\documentclass[10pt,twoside,a4paper,notitlepage]{article}
\usepackage[a4paper,margin=2cm,vmargin={2cm,1.5cm},includeheadfoot]{geometry}

\usepackage{amsmath}
\usepackage{amssymb}
\usepackage{amsfonts}
\usepackage{amsthm}
\usepackage{enumerate}
\usepackage[scaled=0.92]{helvet}
\usepackage{srcltx}
\usepackage{booktabs}
\usepackage{pstricks}

\usepackage{graphicx}
\usepackage[labelfont=bf,margin=1cm,justification=justified,labelsep=period]{caption}
\usepackage{subfigure}

\theoremstyle{plain}
\newtheorem{theorem}{Theorem}[section]
\newtheorem{proposition}[theorem]{Proposition}
\newtheorem{definition}[theorem]{Definition}
\newtheorem{lemma}[theorem]{Lemma}
\newtheorem{corollary}[theorem]{Corollary}

{\theoremstyle{definition}
\newtheorem{remark}[theorem]{Remark}}

\allowdisplaybreaks[1]

\begin{document}

\begin{center}

{\LARGE \bf  L{\'e}vy area for Gaussian processes:\\
 A double Wiener--It\^{o} integral approach} \\ [.2cm]

{\large \bf Albert Ferreiro-Castilla $^a$, Frederic Utzet $^b$
\footnote{Corresponding author. tel +34 935813470: fax: +34 935812790\\
E-mail-addresses: aferreiro@mat.uab.cat (A. Ferreiro-Castilla), utzet@mat.uab.cat (F. Utzet).}}


{\it $^a$ Centre de Recerca Matem{\`a}tica, Campus de Bellaterra, Edifici C,\\
08193 Bellaterra (Barcelona) Spain.}\\
{\it $^b$ Departament de Matem{\`a}tiques, Universitat Aut{\`o}noma de Barcelona, Edifici C,\\   08193 Bellaterra (Barcelona) Spain.}
\end{center}

\noindent\hrule

\bigskip

\noindent{\bf Abstract.}

Let $\{X_{1}(t)\}_{0\leq t\leq1}$ and $\{X_{2}(t)\}_{0\leq t\leq1}$ be two
 independent continuous centered Gaussian processes with  covariance functions
  $R_{1}$ and $R_{2}$. This paper shows that if the covariance functions are of finite
   $p$-variation and $q$-variation respectively and such that $p^{-1}+q^{-1}>1$,
    then the L{\'e}vy area can be defined as a double Wiener--It\`o integral with respect to an isonormal
    Gaussian process induced by  $X_{1}$ and $X_{2}$. Moreover, some properties of the characteristic function
    of that generalised L{\'e}vy area  are studied.

\bigskip


\bigskip


\bigskip

\noindent\hrule

\section{Introduction}

Let $\{W_{1}(t)\ |\ 0\leq t\leq1\}$ and $\{W_{2}(t)\ |\ 0\leq t\leq1\}$ be two independent
standard Wiener processes defined in a probability space $(\Omega,\mathcal{F},P)$,
 and let $A$ be the area included by the curve
\begin{equation*}
x=W_{1}(t)\ ,\quad y=W_{2}(t)\qquad 0\leq t\leq1
\end{equation*}

\noindent and its chord. This random variable was first introduced by L{\'e}vy in \cite{Levy51},
 where it is described by means of stochastic integrals as
\begin{equation*}
\label{definicio}
A=\int_{0}^{1}W_{1}(t)dW_{2}(t)-\int_{0}^{1}W_{2}(t)dW_{1}(t)\ .
\end{equation*}

\noindent L\'evy \cite{Levy51} also computed its characteristic function, which   is
\begin{equation}
\label{characteristic}
\varphi(t):=\mathbb{E}[e^{itA}]=\frac{1}{\cosh(t)}\ ,\qquad t\in\mathbb{R}\ .
\end{equation}

It is easy to show that $A$ has the law of an element of the homogeneous second Wiener chaos generated by
a Brownian motion. This can be proved directly  due to the fact that the elements of the homogeneous second Wiener
chaos have a very particular characteristic function (see Janson \cite[ch. 6]{Jan97}),
and one can realize that (\ref{characteristic}) has this form thanks to  the factorisation
\begin{equation}\label{cosh_decomposition}
\cosh(z)=\prod_{n\geq0}\left(1+\frac{4z^2}{\pi^{2}(2n+1)^{2}}\right)
=\prod_{n\in\mathbb{Z}}\hskip-0.1cm ^{'} \left(1-2iz\alpha_{n}\right)^{1/2}e^{iz\alpha_{n}}\ ,
\end{equation}

\noindent where $\alpha_{n}=(\pi(2n+1))^{-1}$ and $\prod^{'}$ means that each factor is repeated
 twice. Alternatively,  also from (\ref{cosh_decomposition}), the law of  $A$ can be given as the
 law  of a double Wiener--It{\^o} integral
\begin{equation*}
I_{2}^{B}(f):=\iint_{[0,1]^2}f(s,t)dB(s)dB(t)
\end{equation*}

\noindent for an arbitrary Brownian motion $B=\{B_{t}\ |\ t\geq 0\}$, where the kernel
$f$ is obtained in the following way:  consider an orthonormal basis of $L^{2}([0,1])$, which for convenience we write as $\{\phi_{n},\psi_{n}\ |\ n\in\mathbb{Z}\}$, and define
\begin{equation*}
f(s,t)=\sum_{n\in\mathbb{Z}}\alpha_{n}\phi_{n}(s)\phi_{n}(t)+\sum_{n\in\mathbb{Z}}\alpha_{n}\psi_{n}(s)\psi_{n}(t)\ .
\end{equation*}

\noindent Then
\begin{equation*}
I_{2}^{B}(f)=\sum_{n\in\mathbb{Z}}2\alpha_{n}H_{2}(I_{1}^{B}(\phi_{n}))
+\sum_{n\in\mathbb{Z}}2\alpha_{n}H_{2}(I_{1}^{B}(\psi_{n}))\ ,
\end{equation*}

\noindent where $I_{1}^{B}(\phi)=\int_{0}^{1}\phi(s)dB(s)$ and
 $H_{2}(x)=(x^{2}-1)/2$ is the Hermite polynomial of order two.
  Observe that $2 H_{2}(I_{1}^{B}(\phi_{n}))$ and $2 H_{2}(I_{1}^{B}(\psi_{n}))$ are all
  independent centered $\chi^2(1)$ random variables and thus $A\stackrel{\text{law}}{=}I_{2}^{B}(f)$.

However, it is more difficult to get a strong representation of $A$ as a double Wiener--It{\^o}
 integral with respect to the original Brownian motions $W_{1}$ and $W_{2}$. To this end, it
 is necessary to rely on the construction of multiple It{\^o}--Wiener integrals for a general
 white noise, see for instance Nualart \cite[pages  8 and 14]{Nua06}. There, both Brownian
 motions $W_{1}$ and $W_{2}$ are embedded in a Gaussian noise $W$ on $[0,1]\times \{1,2\}$.
  For $h\in L^{2}([0,1]\times\{1,2\},dt\otimes\text{Card})\cong L^{2}([0,1],\mathbb{R}^{2})$,
  where $\text{Card}$ is the counting measure, we have
\begin{equation*}
\int_{[0,1]\times\{1,2\}}hdW=
\int_{0}^{1}h(s,1)dW_{1}(s)+\int_{0}^{1}h(s,2)dW_{2}(s)\ .
\end{equation*}

\noindent Moreover, for $f\in L^{2}(([0,1]\times\{1,2\})^{2})$  symmetric,
\begin{equation}\label{area_chaos}
I_{2}(f)=\sum_{i,j=1}^{2}\iint_{[0,1]^{2}}f((s,i),(t,j))dW_{i}(s)dW_{j}(t)=
2\sum_{i,j=1}^{2}\int_{0}^{1}\int_{0}^{t}f((s,i),(t,j))dW_{i}(s)dW_{j}(t)\ ,
\end{equation}

\noindent see Nualart \cite[p. 23]{Nua06}. For the sake of simplicity, we will indistinctly use $f_{ij}(s,t)$ for $f((s,i),(t,j))$ from now on. Define
\begin{equation}\label{levy_kernel}
f_{ij}^{\mathcal{L}}(s,t)=\left\{
    \begin{array}{ll}
    0,&\text{if } i=j\\
    \frac{1}{2}(\mathbf{1}_{T_{1}}(s,t)-\mathbf{1}_{T_{2}}(s,t)),&\text{if } i=1,j=2\\
    \frac{1}{2}(\mathbf{1}_{T_{2}}(s,t)-\mathbf{1}_{T_{1}}(s,t)),&\text{if } i=2,j=1
    \end{array}
\right.\ ,
\end{equation}

\noindent where $\mathbf{1}_{C}$ is the indicator function of the set $C$ and
\begin{equation*}
T_{1}:=\{(s,t)\in[0,1]^{2}\ |\ s<t\}\ ,\qquad T_{2}:=\{(s,t)\in[0,1]^{2}\ |\ s>t\}\ .
\end{equation*}

\noindent Note that $f^{\mathcal{L}}$ is symmetric, and from (\ref{area_chaos}) and (\ref{definicio})
 it follows that $A\stackrel{\text{a.s.}}{=}I_{2}(f^{\mathcal{L}})$. We will refer to (\ref{levy_kernel}) as the L{\'e}vy kernel.

The aim of this paper is to extend the above strong construction in order to define
 the L{\'e}vy area
 for general Gaussian processes under minimal conditions of their covariance functions,
 and to study its characteristic function.
  We will consider two  independent continuous centered Gaussian processes  $\{X_{1}(t)\ |\ 0\leq t\leq1\}$
  and $\{X_{2}(t)\ |\ 0\leq t\leq1\}$
 with (continuous) covariance functions
  $R_{1}$ and $R_{2}$ and we prove that if the covariance functions are of finite
   $p$-variation and $q$-variation respectively and such that $p^{-1}+q^{-1}>1$,
    then the L{\'e}vy area can be defined as an element of the second Wiener chaos generated
    by $X_{1}$ and $X_{2}$. Such a kind of results have been obtained  (for the non-antisymmetrized
    L{\'e}vy area $\int_0^1X_1(t)\, dX_2(t)$) in the context of rough path analysis by Fritz and
     Victoir, \cite{FV10b,FV10}, but, as far as we know, in such generality they are new for
     classical Gaussian processes.
Our results applied   to two fractional Brownian motions of Hurst parameter
$H$ and $H'$ states that the L{\'e}vy area can be defined if $H+H'>1/2$. In particular, if
$H=H'$, then the condition is $H\in (1/4,1)$ which is a known result (see  Neuenkirch
{\it et al.} \cite{NTU10} and the references
therein) but we present an alternative point of view based in the Huang and Cambanis  \cite{HC78} approach
to stochastic integration for Gaussian processes.
Our results  also extends the ones given by Bardina and Tudor \cite{BT07} where the integral
$\int_0^1X^{H}_t\, dX^{H'}_t$ is defined using Malliavin calculus techniques for
$H\in(0,1)$ and $H'>1/2$.

The paper is organized as follows. We first introduce the general framework of
 the isonormal Gaussian processes, and following the scheme of Huang and Cambanis
 \cite{HC78}, we associate an isonormal Gaussian process to a pair of independent Gaussian processes.
  We also
 give here a definition of a generalised L{\'e}vy area. In the next section we
  derive the conditions on the covariance functions so that $f^{\mathcal{L}}$ generates
   a L{\'e}vy area. As an example we explore what happens with two fractional Brownian
    motion (fBm) with the same covariance function, and the case with two different
     covariance function. This later case allows us to let one of the processes be as
     irregular as desired, that is no low bounds for its Hurst parameter is assumed, provided that the
      other one is regular enough. Finally, we  discuss about the representation of the
      characteristic function of a double Wiener-It\^o integral in terms of a
       Carleman--Fredholm  determinant, that we apply  to compute the characteristic
        function of $A$. Under a further condition of symmetry over the stochastic
         processes we will show that  the characteristic
        function of a  generalised L{\'e}vy area has a Carleman-Fredholm
         determinant with symmetric poles and even multiplicity.

\section{Isonormal Gaussian processes}\label{wiener_chaos_sec}

The framework where (Gaussian) multiple integrals
 are defined is the one of isonormal Gaussian processes. Main reference are Nualart \cite{Nua06},
  and  Peccati and Taqqu \cite{PT10}.
 The more general abstract context of Gaussian Hilbert spaces  developed by Janson \cite{Jan97}
 is also very useful and interesting.

Let $H$ be a separable Hilbert space with inner product $\langle\cdot,\cdot\rangle_{H}$.
An isonormal Gaussian process $\{X(f)\ |\ f\in H\}$ is a centered Gaussian family of random
 variables such that $\mathbb{E}[X(f)X(g)]=\langle f,g\rangle_{H}$. It is well known that the
  construction of the multiple Wiener--It{\^o} integrals with respect to a Brownian motion
  can be transferred to isonormal Gaussian processes; see Nualart \cite[pages 7 and 8]{Nua06}, or  Peccati and Taqqu \cite[ch. 8.4]{PT10}. In that general setup,
    $H^{\otimes n}$ (resp. $H^{\odot n}$) denotes the $n$th (Hilbert) tensor power of $H$
     (resp. the $n$th symmetric tensor power), and $I_{n}(f)$ for $f\in H^{\odot n}$ its $n$th
     multiple integral. For detailed  constructions of that Hilbert spaces see
     Janson \cite{Jan97}.

\subsection{The isonormal Gaussian process associated with  two Gaussian processes}

In this section we describe how  two ordinary
Gaussian process can be imbedded into an isonormal Gaussian process. We extend
  Huang and Cambanis \cite{HC78} approach,  where that construction
was done for one Gaussian process.
Let $X_{1}=\{X_{1}(t)\ |\ t\in[0,1]\}$ and $X_{2}=\{X_{2}(t)\ |\ t\in[0,1]\}$
 be two independent continuous centered Gaussian processes , both starting at zero,  with (continuous)  covariance function
  $R_{1}(s,t)$ and $R_{2}(s,t)$ respectively.
Following Huang and Cambanis \cite{HC78}, let $\mathcal{E}$ denote the set of step functions on $[0,1]$
\begin{equation*}
\phi(t)=\sum_{i=1}^{n}a_{i}\mathbf{1}_{(t_{i},t_{i+1}]}(t)\qquad a_{i}\in\mathbb{R}\ .
\end{equation*}

\noindent Associated with $R_{i}$, for $i=1,2$, we can construct the Hilbert space $H_{i}$ which is the completion of $\mathcal{E}$ under the inner product (with the convenient identifications):
\begin{equation*}
\langle \phi_{1},\phi_{2}\rangle_{H_{i}}:=
\iint_{[0,1]^{2}}\phi_{1}(s)\phi_{2}(t)dR_{i}(s,t)\ .
\end{equation*}

\noindent The above integral is defined so that
\begin{equation}\label{def_int_R}
\iint_{(0,u]\times(0,v]}dR_{i}(s,t)=R_{i}(u,v)\ .
\end{equation}

\begin{remark}\label{cont_cov}
Due to the continuity of the covariance functions the limits of integration in (\ref{def_int_R}) might or might not be included in the integral without changing the result. For instance
\begin{equation*}
\iint_{[0,u]\times[0,v]}dR_{i}(s,t)=\iint_{(0,u]\times(0,v]}dR_{i}(s,t)\ .
\end{equation*}

\end{remark}

In order to define an isonormal Gaussian process associated to both $X_1$ and $X_2$, the set of appropriate elementary functions, $\mathcal{E}_{2}$, are the ones that can be written as $f(t,i)=\delta_{1i}\phi_{1}(t)+\delta_{2i}\phi_{2}(t)$ for $\phi_{1},\ \phi_{2}\in\mathcal{E}$, where $\delta_{ij}$ is th Kronecker's delta. It is clear that on $\mathcal{E}_{2}$ we can consider the inner product (with the convenient identifications):

$$\langle f,g\rangle_{\mathcal{E}_{2}}
=\langle f(\cdot,1),g(\cdot,1)\rangle_{H_{1}}+\langle f(\cdot,2),g(\cdot,2)\rangle_{H_{2}}
=\iint_{[0,1]^{2}}f(s,1)g(t,1)dR_{1}(s,t)+\iint_{[0,1]^{2}}f(s,2)g(t,2)dR_{2}(s,t)\ .$$

\noindent Let us call $H$ the Hilbert space which is the completion of $\mathcal{E}_{2}$ with the above inner product. Next lemma characterises $H$; its proof is straightforward.

\begin{lemma}\label{lemma_iso}
Under the above notation $H\cong H_{1}\oplus H_{2}$, where $H_{1}\oplus H_{2}$ is the
Hilbertian direct sum
of $H_1$ and $H_2$, that is the Hilbert space which consists in all ordered pairs $(x_{1},x_{2})\in H_{1}\times H_{2}$ equipped with the inner product
$
\langle(x_{1},x_{2}),(y_{1},y_{2})\rangle_{H_{1}\oplus H_{2}}=
\langle x_{1},y_{1}\rangle_{H_{1}}+
\langle x_{2},y_{2}\rangle_{H_{2}}
$.

\end{lemma}

\noindent Now we are ready to construct the isonormal Gaussian process indexed by $H$
 which will represent the $2$-dimensional process $\{(X_{1}(t),X_{2}(t))\ |\ t\in[0,1]\}$.
  From the independence of $X_{1}$ and $X_{2}$, it turns out that
  $X:\mathcal{E}_{2}\to L^{2}(\Omega,\mathcal{F},P)$ defined by
$
X(f):=X_{1}(f(\cdot,1))+X_{2}(f(\cdot,2))
$
is an isometry which can be extended to $H$. Thus $X=\{X(f)\ |\ f\in H\}$ is an
 isonormal Gaussian process.
\subsection{Generalised L{\'e}vy area}
 In the previous  context,  the generalised L{\'e}vy area will be an element of the second Wiener chaos with respect to the process $X$. Therefore we need to identify the elements of $H^{\otimes 2}$. Note that
\begin{equation}\label{factorisation_H}
H^{\otimes 2}\cong H^{\otimes 2}_{1}\oplus (H_{1}\otimes H_{2})\oplus
(H_{2}\otimes H_{1})\oplus H^{\otimes 2}_{2}\ ,
\end{equation}

\noindent which gives a very appropriate interpretation of the elements $f\in H^{\otimes 2}$ as $2$ by $2$ matrices with entries $f_{ij}\in H_{i}\otimes H_{j}$ for $i,j=1,2$. The above isometry also induces the decomposition
\begin{equation*}
\langle f,g\rangle_{H^{\otimes 2}}=
\sum_{i,j=1}^{2}\langle f_{ij},g_{ij}\rangle_{H_{i}\otimes H_{j}}\ .
\end{equation*}

The double It{\^o}--Wiener integral, $I_{2}(\cdot)$, is an isometry between $H^{\odot 2}$ and
 the second Wiener chaos. Therefore, the desirable definition of a generalised L{\'e}vy area
  would be $I_{2}(f^{\mathcal{L}})$ whenever $f^{\mathcal{L}}\in H^{\odot 2}$,
  where $f^{\mathcal{L}}$ was defined in (\ref{levy_kernel}). Unfortunately this is very difficult,
   if possible at all, to prove. We will circumvent this problem by finding an element in
   $H^{\odot 2}$ which is indistinguishable from $f^{\mathcal{L}}$ and to which we will apply
   the isometry $I_{2}(\cdot)$.
 In other words, we will say that a function $f\in L^{2}(([0,1]\times\{1,2\})^{2})$ (symmetric)
  belongs to $H^{\odot 2}$ as long as there is an element $\hat{f}\in H^{\odot 2}$ such that
\begin{equation*}
\langle\hat{f},g\rangle_{H^{\otimes 2}}=
\sum_{i,j=1}^{2}\iiiint_{[0,1]^{4}}f_{ij}(s,t)g_{ij}(u,v)dR_{i}(s,u)dR_{j}(t,v)\qquad \forall g\in\mathcal{E}_{2}\ ,
\end{equation*}

\noindent where the above integral is an iterated Riemann--Stieltjes integral,
see Corollary \ref{approx_norm} below. Note that we are not enlarging the space
$H^{\odot 2}$ but renaming the element $\hat{f}$ by $f$, since under the inner product in
$H^{\odot 2}$ they are indistinguishable. Thus the map $I_{2}$ is well defined for $\hat{f}$
 and we put $I_{2}(f):=I_{2}(\hat{f})$. This is a common procedure to ease the identification
 of the elements of Hilbert spaces which have been constructed by completion, see Huang and
  Cambanis \cite{HC78}. Therefore we will define the generalised L{\'e}vy area  in the following
  way:

\begin{definition}[Generalised L{\'e}vy area]\label{GLA}
We will say that $I_{2}(f^{\mathcal{L}})$ is a generalised L{\'e}vy area if there exists
 $\hat{f}^{\mathcal{L}}\in H^{\odot 2}$ such that
\begin{equation}\label{def_area}
\langle\hat{f}^{\mathcal{L}},g\rangle_{H^{\otimes 2}}=
\sum_{i,j=1}^{2}\iiiint_{[0,1]^{4}}f_{ij}^{\mathcal{L}}(s,u)g_{ij}(t,v)dR_{i}(s,t)dR_{j}(u,v)\ ,
\end{equation}

\noindent for all step functions $g\in \mathcal{E}_{2}^{\otimes 2}$. Then $I_{2}(f^{\mathcal{L}}):=I_{2}(\hat{f}^{\mathcal{L}})$.
\end{definition}

Another problem we have to face in order to make this definition tractable is that
we do not know how to compute the inner product
 $\langle\hat{f}^{\mathcal{L}},g\rangle_{H^{\otimes 2}}$.
  We were very careful to only write the integral form of the inner product
   in $H^{\otimes 2}$ for step functions, indeed we only know how to calculate
   the inner product for step functions since the rest of the space was
   constructed by completion. Hence we will need to approximate
   $\hat{f}^{\mathcal{L}}$ by step functions and check equation (\ref{def_area}) as a limit
   equality. Before that, let us explicit the inner product in $H^{\otimes 2}$ for
    step functions as integrals with respect to the covariance functions $R_{1}$ and $R_{2}$.

\begin{lemma}
Let $f,g\in\mathcal{E}^{\otimes 2}$, then
\begin{equation*}
\langle f,g\rangle_{H_{i}\otimes H_{j}}
=\iiiint_{[0,1]^{4}}f(s,t)g(u,v)dR_{i}(s,u)dR_{j}(t,v)\ .
\end{equation*}

\end{lemma}

\noindent{\it Proof.}
Let $f=f_{1}\otimes f_{2}$ and $g=g_{1}\otimes g_{2}$, then
\begin{eqnarray*}
\langle f,g\rangle_{H_{i}\otimes H_{j}}
&=&\langle f_{1},g_{1}\rangle_{H_{i}}\langle f_{2},g_{2}\rangle_{H_{j}}
=\iint_{[0,1]^{2}}f_{1}(s)g_{1}(t)dR_{1}(s,t)\iint_{[0,1]^{2}}f_{2}(u)g_{2}(v)dR_{2}(u,v)\\
&=&\iiiint_{[0,1]^{4}}f_{1}(s)g_{1}(t)f_{2}(u)g_{2}(v)dR_{i}(s,t)dR_{j}(u,v)\ ,
\end{eqnarray*}

\noindent and we get the result since a realization of the tensor product for functions is just the plain product. \quad $\blacksquare$

\begin{corollary}\label{approx_norm}
Let $f,g\in H^{\otimes 2}$, such that $f_{ij},g_{ij}\in\mathcal{E}^{\otimes 2}$. Then
\begin{eqnarray*}
\langle f,g\rangle_{H^{\otimes 2}}
&=&\sum_{i,j=1}^{2}\iiiint_{[0,1]^{4}}f_{ij}(s,u)g_{ij}(t,v)dR_{i}(s,t)dR_{i}(u,v)\\
&=&\sum_{i,j=1}^{2}\iint_{[0,1]^{2}}\left(\iint_{[0,1]^{2}}f_{ij}(s,u)g_{ij}(t,v)dR_{i}(s,t)\right)dR_{i}(u,v)\ .
\end{eqnarray*}

\end{corollary}

\section{Existence of a generalised L{\'e}vy area}

This section will give the sufficient conditions on the processes $\{X_{1}(t)\ |\ 0\leq t\leq1\}$ and $\{X_{2}(t)\ |\ 0\leq t\leq1\}$ so a generalised L{\'e}vy area exists (see Definition \ref{GLA}). In fact the conditions on the processes will be constrains on their covariance functions, indeed this is what the previous sections suggest as the Hilbert space of the domain of $I_{2}(\cdot)$ is characterised by the covariance function of $X_{1}$ and $X_{2}$. Friz and Victoir \cite{FV10} claim the $p$-variation of the covariance function of a Gaussian process to be a fundamental quantity related to the process. Therefore, we first recall some definitions on the $p$-variation of a function.

\subsection{Functions of finite $p$-variation and Young's inequality}

For the sake of completeness and to introduce notation, we here give some definitions on the $p$-variation of a function, what is meant by a control map and state Young's inequality. For further reading see Dudley and Norvai{\v{s}}a \cite{DN99} and Friz and Victoir \cite{FV10,FV10b}.

For a given interval of the real line $[s,t]$ such that $s\leq t$, we will denote the set of all partitions of $[s,t]$ by
\begin{equation*}
\mathcal{P}([s,t]):=\{
\{t_{0},\ldots,t_{n}\}\ |\ s=t_{0}<t_{1}<\ldots<t_{n}=t,\ n\in\mathbb{N}\ \}\ .
\end{equation*}

\noindent If $D\in\mathcal{P}([s,t])$, then write $|D|:=\max_{t_{i}\in D}\{|t_{i}-t_{i-1}|\}$.
\begin{definition}
Let $f:[s,t]\to\mathbb{R}$ be a function and $p\geq 1$. We say that $f$ has finite $p$-variation if $V^{1}_{p}(f,[s,t])<\infty$, where
\begin{equation*}
V^{1}_{p}(f,[s,t]):=\sup_{D\in\mathcal{P}([s,t])}\left(\sum_{t_{i}\in D}|f(t_{i+1})-f(t_{i})|^{p}\right)^{1/p}\ .
\end{equation*}

\end{definition}

The superscript on $V^{1}_{p}$ is to emphasise that $f$ is $1$-dimensional in contrast to the $p$-variation of a $2$-dimensional function, which it is defined below.

\begin{definition}
Let $f:[s,t]\times[u,v]\to\mathbb{R}$ be a function and $p\geq 1$. We say that $f$ has finite $p$-variation if $V^{2}_{p}(f,[s,t]\times[u,v])<\infty$, where
\begin{equation*}
V^{2}_{p}(f,[s,t]\times[u,v]):=
\sup_{\begin{smallmatrix}D\in\mathcal{P}([s,t])\\D'\in\mathcal{P}([u,v])\end{smallmatrix}}
\left(\sum_{(t_{i},t_{j}')\in D\times D'}
\left|f\left(\begin{array}{c}t_{i}\\t_{i+1}\end{array},\begin{array}{c}t_{j}'\\t_{j+1}'\end{array}\right)\right|^{p}\right)^{1/p}
\end{equation*}

\noindent and
\begin{equation*}
f\left(\begin{array}{c}t_{i}\\t_{i+1}\end{array},\begin{array}{c}t_{j}'\\t_{j+1}'\end{array}\right):=
f(t_{i+1},t_{j+1}')-f(t_{i+1},t_{j}')-f(t_{i},t_{j+1}')+f(t_{i},t_{j}')\ .
\end{equation*}

\end{definition}

Another important concept related to the $p$-variation is the \emph{control map} (see Friz and Victoir \cite[ch. 5.1, 5.3]{FV10b}).

\begin{definition}\label{def_control}
A $2$-dimensional control is a map $\omega$ from $[s,t]\times[u,v]$ to $[0,\infty)$ where $0\leq s\leq t\leq 1$, $0\leq u\leq v\leq1$ and such that for all $r\leq s\leq t$, $u\leq v$,
\begin{eqnarray*}
\omega ([r,s]\times[u,v])+\omega ([s,t]\times[u,v])&\leq&\omega([r,t]\times[u,v])\\
\omega ([u,v]\times[r,s])+\omega ([u,v]\times[s,t])&\leq&\omega([u,v]\times[r,t])\ ,
\end{eqnarray*}

\noindent and $\lim_{s\to t}\omega([s,t]\times[u,v])=\lim_{s\to t}\omega([u,v]\times[s,t])=0$.

\end{definition}

It is just for convenience that we set the variables to be in $[0,1]$. The relationship between the control and the $p$-variation is given by the following lemma (Friz and Victoir \cite[p. 106]{FV10b}):

\begin{lemma}\label{variation_control}
Let $f$ be a continuous function of finite $p$-variation -- $V^{2}_{p}(f,[s,t]\times[u,v])<\infty$ --, then there is a $2$--dimensional control map, $\omega$, such that
\begin{equation*}
V^{2}_{p}(f,[s,t]\times[u,v])\leq \omega^{1/p}([s,t]\times[u,v])\ .
\end{equation*}

\end{lemma}

We  will need the following technical result about
 the product of control maps; see Fritz and Victoir \cite{FV10};  its proof is a consequence
 of a discrete H{\"o}lder type inequality proved by Young \cite[p. 252]{Young36}.

\begin{lemma}\label{lemma_control_prod}
Let $\omega_{1}$ and $\omega_{2}$ be $2$-dimensional control maps over the same rectangle and $p,q>0$ such that $p^{-1}+q^{-1}\geq1$, then $\omega_{1}^{1/p}\omega_{2}^{1/q}$ is also a $2$-dimensional control map.
\end{lemma}

Finally, we recall the statement of Young's inequality for a $2$-dimensional function (see Towghi \cite{Towghi02}):

\begin{theorem}\label{YS_integral}
Let $f$ and $g$ be functions such that
\begin{enumerate}[i)]
\item $V^{2}_{p}(f,[0,1]\times[0,1])<\infty$,
\item $V^{1}_{p}(f(0,\cdot),[0,1])<\infty$,
\item $V^{1}_{p}(f(\cdot,0),[0,1])<\infty$,
\item $|f(0,0)|<\infty$,
\item $V^{2}_{p}(g,[0,1]\times[0,1])<\infty$
\end{enumerate}

\noindent and $p^{-1}+q^{-1}>1$, and define
\begin{equation*}
||f||_{W_{p}^{2}([0,1]^{2})}:=
V^{2}_{p}(f,[0,1]\times[0,1])+V^{1}_{p}(f(0,\cdot),[0,1])+V^{1}_{p}(f(\cdot,0),[0,1])+|f(0,0)|\ .
\end{equation*}

\noindent If $f$ and $g$ do not have any common jump points then the Young--Stieltjes integral of $f$ with respect to $g$ exists, and
\begin{equation*}
\left|\int_{0}^{1}\int_{0}^{1}f(x,y)dg(x,y)\right|\leq
c(p,q)||f||_{W_{p}^{2}([0,1]^{2})}V^{2}_{p}(g,[0,1]\times[0,1])\ ,
\end{equation*}

\noindent where $c(p,q)$ is a constant independent of $f$ and $g$.

\end{theorem}

In our setting, the functions to which we are going to apply Young's inequality are continuous and thus do not have jump points.

\begin{remark}
The definition of finite $p$--variation could be stated for $p>0$ both in the
 $1$--dimensional and in the $2$--dimensional case, but we restrict ourselves
  to $p\geq 1$. This is because a $1$--dimensional continuous function of finite $p$--variation
  for $p<1$ is constant (see Friz and Victoir \cite[p. 78]{FV10}).
   This is not true for the $2$--dimensional case, for example the function $f(x,y)=x+y$
    has finite $p$--variation for all $p>0$. However, for continuous covariance functions
    coming from processes that start at a point rather than from a distribution
    it is true (see next result).
    We will see in the next section that the hypotheses related with
    the finite variation are always with respect to continuous covariance functions.
     Therefore, without lost of generality, we consider
      $p$--variations for $p\geq 1$.
\end{remark}

\begin{lemma}
A continuous function $f$ on $[0,1]^{2}$ such that $f(0,0)=f(s,0)=f(0,t)$ for all $t,s\in[0,1]$ and of finite $p$--variation with $p<1$ is constant.
\end{lemma}

\noindent{\it Proof.}
For a fixed $a\in[0,1]$,  the function $y\to f(a,y)-f(0,y)$ is a $1$--dimensional continuous
function of finite $p$--variation and hence constant. Indeed, it is zero since $f(a,0)=f(0,0)$,
and the result follows.
\quad $\blacksquare$

\subsection{Main result}

The main result of the paper is proved in this section. We construct a sequence of step
 functions which converge almost sure to the L{\'e}vy kernel and show that it is a Cauchy s
 equence in $H^{\otimes 2}$, and finally, we also show that its limit satisfies Definition
 \ref{GLA}.

We start by a technical lemma which will ease the proof of the main result, but, before that,
 let us introduce some notation which will be used extensively in this section. According to
 definition (\ref{def_int_R}) we have that
\begin{equation*}
\int_{s}^{t}\int_{u}^{v}dR_{i}(x,y)=
R_{i}\left(\begin{array}{c}s\\t\end{array},\begin{array}{c}u\\v\end{array}\right)\ .
\end{equation*}

\noindent Let $\{t_{i}^{n}=i2^{-n}\ |\ i=0,\ldots,2^{n}-1\}$ be the dyadic partition of the interval $[0,1]$ for a given $n$, and consider the dyadic partition of the triangles $T_{1}$ and $T_{2}$
\begin{equation*}
T_{1}^{n}:=\bigcup_{i<j}I_{i}^{n}\times I_{j}^{n}\qquad
T_{2}^{n}:=\bigcup_{i>j}I_{i}^{n}\times I_{j}^{n}\ ,
\end{equation*}

\noindent where $I_{i}^{n}:=(t_{i}^{n},t_{i+1}^{n}]$. Then a natural approximation of the
 L{\'e}vy kernel will be
\begin{equation*}
f_{n}((s,i),(t,j)):=\left\{
    \begin{array}{ll}
    0,&\text{if } i=j\\
    \frac{1}{2}(\mathbf{1}_{T_{1}^{n}}(s,t)-\mathbf{1}_{T_{2}^{n}}(s,t)),&\text{if } i=1,j=2\\
    \frac{1}{2}(\mathbf{1}_{T_{2}^{n}}(s,t)-\mathbf{1}_{T_{1}^{n}}(s,t)),&\text{if } i=2,j=1
    \end{array}
\right.\ .
\end{equation*}

\begin{lemma}\label{cauchy_technical}
Let $R_{1}$ and $R_{2}$ be two continuous covariance functions in $[0,1]^{2}$. Let $R_{1}$ be of finite $p$-variation and $R_{2}$ of finite $q$-variation and assume that \mbox{$p^{-1}+q^{-1}>1$}, then
\begin{equation*}
\lim_{\begin{smallmatrix}n\to\infty\\m\to\infty\end{smallmatrix}}
\underbrace{\sum_{i,j=1}^{2}\iiiint_{[0,1]^{4}}(f^{\mathcal{L}}-f_{n})_{ij}(s,u) \cdot (f^{\mathcal{L}}-f_{m})_{ij}(t,v)dR_{i}(s,t)dR_{j}(u,v)}_{(\star)}=0\ .
\end{equation*}

\end{lemma}

\noindent{\it Proof.}
Write $J_{k,l}^{n,m}:=I_{k}^{n}\times I_{l}^{m}$ and note that
\begin{equation*}
(f^{\mathcal{L}}-f_{n})_{ij}(s,u) \cdot (f^{\mathcal{L}}-f_{m})_{ij}(t,v)=\frac{(1-\delta_{ij})}{4}
\sum_{k=0}^{2^{n}-1}\sum_{l=0}^{2^{m}-1}\mathbf{1}_{(J_{k,l}^{n,m})^{2}}(s,t,u,v)
    (\mathbf{1}_{\{(v-t)(u-s)>0\}}-\mathbf{1}_{\{(v-t)(u-s)<0\}})\ .
\end{equation*}

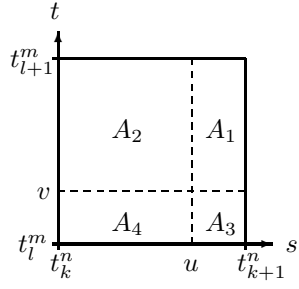
\begin{figure}
\centering
\begin{picture}(100,90)(0,0)
\linethickness{0.3mm}
\put(18,10){\line(1,0){82}}
\put(100,10){\vector(1,0){0.12}}
\linethickness{0.3mm}
\put(20,8){\line(0,1){82}}
\put(105,7){$s$}
\put(17,95){$t$}
\put(67,0){$u$}
\put(87,0){$t_{k+1}^{n}$}
\put(17,0){$t_{k}^{n}$}
\put(12,27){$v$}
\put(5,7){$t_{l}^{m}$}
\put(2,77){$t_{l+1}^{m}$}
\put(40,50){$A_{2}$}
\put(75,50){$A_{1}$}
\put(40,15){$A_{4}$}
\put(75,15){$A_{3}$}
\put(20,90){\vector(0,1){0.12}}
\linethickness{0.3mm}
\put(18,80){\line(1,0){72}}
\linethickness{0.3mm}
\put(90,8){\line(0,1){72}}
\linethickness{0.1mm}
\multiput(70,10)(0,4.83){15}{\line(0,1){2.41}}
\linethickness{0.1mm}
\multiput(20,30)(4.83,0){15}{\line(1,0){2.41}}
\end{picture}
\caption{Integration of the $(\star\star)$-term over the single rectangle $J_{k,l}^{n,m}$ for a fixed $(u,v)\in J_{k,l}^{n,m}$.\label{esquema}}
\end{figure}

\noindent Therefore the quadruple integral of the $(\star)$-term is split into a sum
of quadruple integrals over $(J_{k,l}^{n,m})^2$. These integrals are iterated
 integrals and they can be further reduced, according to Figure \ref{esquema}, to
\begin{equation*}
(\star)=\frac{1}{4}\sum_{\begin{smallmatrix}i,j=1\\i\neq j\end{smallmatrix}}^{2}
    \underbrace{\sum_{k=0}^{2^{n}-1}\sum_{l=0}^{2^{m}-1}
    \iint_{J_{k,l}^{n,m}}F_{k,l}^{n,m,i}(u,v)dR_{j}(u,v)}_{(\star\star)}\ ,
\end{equation*}

\noindent where
\begin{eqnarray*}
F_{k,l}^{n,m,i}(u,v)&:=&
\iint_{A_{1}}dR_{i}(s,t)-\iint_{A_{2}}dR_{i}(s,t)-\iint_{A_{3}}dR_{i}(s,t)+\iint_{A_{4}}dR_{i}(s,t)\\
&=&R_{i}\left(\begin{array}{c}u\\t_{k+1}^{n}\end{array},\begin{array}{c}v\\t_{l+1}^{m}\end{array}\right)-
R_{i}\left(\begin{array}{c}t_{k}^{n}\\u\end{array},\begin{array}{c}v\\t_{l+1}^{m}\end{array}\right)-
R_{i}\left(\begin{array}{c}u\\t_{k+1}^{n}\end{array},\begin{array}{c}t_{l}^{m}\\v\end{array}\right)+
R_{i}\left(\begin{array}{c}t_{k}^{n}\\u\end{array},\begin{array}{c}t_{l}^{m}\\v\end{array}\right)\ .
\end{eqnarray*}

\noindent Note that in the above definition we have used Remark \ref{cont_cov}.

It is enough to prove that the $(\star\star)$-term goes to zero as $n,m\to\infty$ for $i=1$
 and $j=2$. The key point is to apply Young's inequality to each summand of the $(\star\star)$-term, in that way we will prove the existence of the integrals and get a bound for them. In order to do so note the following identities which relate the function $F_{k,l}^{n,m,1}$ with the function $R_{1}$
\begin{equation*}
\begin{array}{rclcrcl}
F_{k,l}^{n,m,1}\left(\begin{array}{c}u\\u'\end{array},\begin{array}{c}v\\v'\end{array}\right)
    &=&4R_{1}\left(\begin{array}{c}u\\u'\end{array},\begin{array}{c}v\\v'\end{array}\right)
&,&F_{k,l}^{n,m,1}(t_{k}^{n},v)-F_{k,l}^{n,m,1}(t_{k}^{n},v')
    &=&2R_{1}\left(\begin{array}{c}t_{k}^{n}\\t_{k+1}^{n}\end{array},\begin{array}{c}v\\v'\end{array}\right)\ ,\\
\\
F_{k,l}^{n,m,1}(t_{k}^{n},t_{l}^{m})
    &=&R_{1}\left(\begin{array}{c}t_{k}^{n}\\t_{k+1}^{n}\end{array},\begin{array}{c}t_{l}^{m}\\t_{l+1}^{m}\end{array}\right)
&\text{ and }&F_{k,l}^{n,m,1}(u,t_{l}^{m})-F_{k,l}^{n,m,1}(u',t_{l}^{m})
    &=&2R_{1}\left(\begin{array}{c}u\\u'\end{array},\begin{array}{c}t_{l}^{m}\\t_{l+1}^{m}\end{array}\right)\ .
\end{array}
\end{equation*}

\noindent Now, since $p^{-1}+q^{-1}>1$, there exists $r,p'>0$ such that $r+(p')^{-1}=p^{-1}$ and
 $(p')^{-1}+q^{-1}>1$. Thus, from $p'>p\geq 1$ we have that $R_{1}$ is also of finite $p'$--variation and use the above equalities to obtain
\begin{eqnarray*}
\begin{array}{rclcrcl}
V_{p'}^{2}(F_{k,l}^{n,m,1},\overline{J_{k,l}^{n,m}})
    &=&4V_{p'}^{2}(R_{1},\overline{J_{k,l}^{n,m}})
&,&V_{p'}^{1}(F_{k,l}^{n,m,1}(t_{k}^{n},\cdot),[t_{l}^{m},t_{l+1}^{m}])
    &\leq&2V_{p'}^{2}(R_{1},\overline{J_{k,l}^{n,m}})\ ,\\
\left|F_{k,l}^{n,m,1}(t_{k}^{n},t_{l}^{m})\right|
    &\leq&V_{p'}^{2}(R_{1},\overline{J_{k,l}^{n,m}})
&\text{ and }&V_{p'}^{1}(F_{k,l}^{n,m,1}(\cdot,t_{l}^{m}),[t_{k}^{n},t_{k+1}^{n}])
    &\leq&2V_{p'}^{2}(R_{1},\overline{J_{k,l}^{n,m}})\ ,
\end{array}
\end{eqnarray*}

\noindent where $\overline{J_{k,l}^{n,m}}$ is the closure of $J_{k,l}^{n,m}$.
This consideration is a technicality required by the definition of the finite variation.
According to the notation of Theorem \ref{YS_integral}
 the above inequalities imply that $||F_{k,l}^{n,m,1}||_{W_{p'}^{2}(\overline{J_{k,l}^{n,m}})}\leq 9V_{p'}^{2}(R_{1},\overline{J_{k,l}^{n,m}})$, and thus we can apply Young's inequality to every integral of the $(\star\star)$-term to get
$$
|(\star\star)|
\leq c(p',q)\sum_{k=0}^{2^{n}-1}\sum_{l=0}^{2^{m}-1}
||F_{k,l}^{n,m,1}||_{W_{p'}^{2}(\overline{J_{k,l}^{n,m}})}V^{2}_{q}(R_{2},\overline{J_{k,l}^{n,m}})
\leq 9c(p',q)\sum_{k=0}^{2^{n}-1}\sum_{l=0}^{2^{m}-1}
V_{p'}^{2}(R_{1},\overline{J_{k,l}^{n,m}})V^{2}_{q}(R_{2},\overline{J_{k,l}^{n,m}})\ .
$$


Recall Lemma \ref{variation_control} to associate the finite $p$--variation of $R_{1}$ to a control map $\omega_{1}$, and denote by $\omega_{2}$ the control map with respect to the $q$--variation of $R_{2}$. Finally we use Lemma \ref{lemma_control_prod} to deduce that $\hat{\omega}:=\omega_{1}^{1/p'}\omega_{1}^{1/q}$ is a control map which bounds the $(\star\star)$-term as
\begin{eqnarray*}
|(\star\star)|&\leq&C\sum_{k=0}^{2^{n}-1}\sum_{l=0}^{2^{m}-1}
\left(
\sup_{\begin{smallmatrix}D\in\mathcal{P}(\overline{I_{k}^{n}})\\D'\in\mathcal{P}(\overline{I_{l}^{m}})\end{smallmatrix}}
\left(\sum_{(t_{i},t_{j}')\in D\times D'}
\left|R_{1}\left(\begin{array}{c}t_{i}\\t_{i+1}\end{array},\begin{array}{c}t_{j}'\\t_{j+1}'\end{array}\right)\right|^{p'}\right)^{1/p'}
\right)
V^{2}_{q}(R_{2},\overline{J_{k,l}^{n,m}})\\
&\leq&C\sup_{\begin{smallmatrix}|u-u'|\leq 2^{-n}\\|v-v'|\leq 2^{-m}\end{smallmatrix}}
\left|R_{1}\left(\begin{array}{c}u\\u'\end{array},\begin{array}{c}v\\v'\end{array}\right)\right|^{\frac{p'-p}{p'}}
\sum_{k=0}^{2^{n}-1}\sum_{l=0}^{2^{m}-1}
(V_{p}^{2}(R_{1},\overline{J_{k,l}^{n,m}}))^{p/p'}V^{2}_{q}(R_{2},\overline{J_{k,l}^{n,m}})\\
&\leq&C\sup_{\begin{smallmatrix}|u-u'|\leq 2^{-n}\\|v-v'|\leq 2^{-m}\end{smallmatrix}}
\left|R_{1}\left(\begin{array}{c}u\\u'\end{array},\begin{array}{c}v\\v'\end{array}\right)\right|^{\frac{p'-p}{p'}}
\sum_{k=0}^{2^{n}-1}\sum_{l=0}^{2^{m}-1}
\omega_{1}^{1/p'}(\overline{J_{k,l}^{n,m}})\omega_{2}^{1/q}(\overline{J_{k,l}^{n,m}})\\
&\leq&C\sup_{\begin{smallmatrix}|u-u'|\leq 2^{-n}\\|v-v'|\leq 2^{-m}\end{smallmatrix}}
\left|R_{1}\left(\begin{array}{c}u\\u'\end{array},\begin{array}{c}v\\v'\end{array}\right)\right|^{\frac{p'-p}{p'}}
\sum_{k=0}^{2^{n}-1}\sum_{l=0}^{2^{m}-1}
\hat{\omega}(\overline{J_{k,l}^{n,m}})\\
&\leq&C\hat{\omega}([0,1]^{2})\sup_{\begin{smallmatrix}|u-u'|\leq 2^{-n}\\|v-v'|\leq 2^{-m}\end{smallmatrix}}
\left|R_{1}\left(\begin{array}{c}u\\u'\end{array},\begin{array}{c}v\\v'\end{array}\right)\right|^{\frac{p'-p}{p'}}\ ,
\end{eqnarray*}

\noindent where $C$ is a constant which is renamed when necessary. Finally note
that the last supremum goes to zero as $n,m\to\infty$ by the uniform continuity
 of $R_{1}$ and the result follows.\quad $\blacksquare$

Now we are ready to finalize the construction of the generalized L{\'e}vy area,
 we need to prove that the sequence $f_{n}$ in Lemma \ref{cauchy_technical}
 is a Cauchy sequence and its limit, $\hat{f}^{\mathcal{L}}$, satisfies Definition \ref{GLA}.

\begin{theorem}\label{cauchy_aproximation}
Let $\{X_{1}(t)\ |\ 0\leq t \leq1\}$ and $\{X_{2}(t)\ |\ 0\leq t \leq1\}$ be
 two continuous centered Gaussian processes, such that $X_1(0)=X_2(0)=0$, independent,
   and with  covariance functions $R_{1}$ and $R_{2}$
 respectively. Let $R_{1}$ be of finite $p$-variation and $R_{2}$ be of finite
 $q$-variation and assume that $p^{-1}+q^{-1}>1$, then the sequence $\{f_{n}\}_{n\ge 1}$ is a
  Cauchy sequence in $H^{\otimes 2}$. We will denote its limit by $\hat{f}^{\mathcal{L}}$. Moreover, we have that
\begin{equation}\label{def_theorem}
\langle \hat{f}^{\mathcal{L}},g\rangle_{H^{\otimes 2}}
=\sum_{i,j=1}^{2}\iiiint_{[0,1]^{4}}f_{ij}^{\mathcal{L}}(s,u)g_{ij}(t,v)dR_{i}(s,t)dR_{j}(u,v)
\end{equation}

\noindent for all step functions $g\in H^{\otimes 2}$.

\end{theorem}

\noindent{\it Proof.}
Note that $f_{n}-f_{m}\in H^{\otimes 2}$ is a difference of two step functions
 and hence it is a step function itself. Therefore by Lemma \ref{approx_norm} we have that
\begin{eqnarray*}
||f_{n}-f_{m}||_{H^{\otimes 2}}
&=&\sum_{\begin{smallmatrix}i,j=1\\i\neq j\end{smallmatrix}}^{2}
    \iiiint_{[0,1]^{4}}(f_{n}-f_{m})_{ij}(s,u)
        \cdot(f_{n}-f_{m})_{ij}(t,v)dR_{i}(s,t)dR_{i}(u,v)\\
&=&\sum_{\begin{smallmatrix}i,j=1\\i\neq j\end{smallmatrix}}^{2}
    \iiiint_{[0,1]^{4}}(f_{n}-f^{\mathcal{L}}+f^{\mathcal{L}}-f_{m})_{ij}(s,u)
    \cdot(f_{n}-f^{\mathcal{L}}+f^{\mathcal{L}}-f_{m})_{ij}(t,v)dR_{i}(s,t)dR_{i}(u,v)\ .
\end{eqnarray*}

\noindent Each term of the above product was denoted as a $(\star)$-term in Lemma \ref{cauchy_technical} and thus goes to zero as $n,m\to\infty$.

For the second part of the proof it suffices to prove the equality for a function $g$ such that
$g_{12}(s,t)=\mathbf{1}_{[a,b]\times[c,d]}(s,t)$ where $[a,b]\times[c,d]\subseteq[0,1]^{2}$
and $g_{ij}(s,t)\equiv 0$ for $i\neq 1$ or $j\neq 2$. Since $\hat{f}^{\mathcal{L}}$
 is the limit of $\{f_{n}\}_{n\ge 1}$ in $H^{\otimes 2}$ then
\begin{equation*}
\lim_{n\to\infty}\langle f_{n},g\rangle_{H^{\otimes 2}}=
\langle \hat{f}^{\mathcal{L}},g\rangle_{H^{\otimes 2}}\ .
\end{equation*}

\noindent Our objective is to prove that $\lim_{n\to\infty}\langle f_{n},g\rangle_{H^{\otimes 2}}$
equals the left hand side of equation (\ref{def_theorem}). From the definition of
the L{\'e}vy kernel we have that
\begin{eqnarray}
\sum_{i,j=1}^{2}\iiiint_{[0,1]^{4}}f_{ij}^{\mathcal{L}}(s,u)g_{ij}(t,v)dR_{i}(s,t)dR_{j}(u,v)\nonumber\\
&&\hskip-3cm=\frac{1}{2}\iiiint_{[0,1]^{4}}(\mathbf{1}_{s>u}(s,u)-\mathbf{1}_{s<u}(s,u))\mathbf{1}_{[a,b]\times[c,d]}(t,v)dR_{1}(s,t)dR_{2}(u,v)\nonumber\\
&&\hskip-3cm=\frac{1}{2}\iint_{[0,1]\times[c,d]}dR_{2}(u,v)\iint_{[0,1]\times[a,b]}(\mathbf{1}_{s>u}(s,u)-\mathbf{1}_{s<u}(s,u))dR_{1}(s,t)\nonumber\\
&&\hskip-3cm=\frac{1}{2}\iint_{[0,1]\times[c,d]}
\left[R_{1}\left(\begin{array}{c}u\\1\end{array},\begin{array}{c}a\\b\end{array}\right)-
R_{1}\left(\begin{array}{c}0\\u\end{array},\begin{array}{c}a\\b\end{array}\right)\right]
dR_{2}(u,v)\label{int_aux_1}\ .
\end{eqnarray}

\noindent The above integral  is a well defined  Young--Stieltjes integral. Then, for
 $D\in\mathcal{P}([0,1])$ and $D'\in\mathcal{P}([c,d])$,
\begin{eqnarray}
\sum_{i,j=1}^{2}\iiiint_{[0,1]^{4}}f_{ij}^{\mathcal{L}}(s,u)g_{ij}(t,v)dR_{i}(s,t)dR_{j}(u,v)\nonumber\\
&&\hskip-4cm=\frac{1}{2}
\lim_{\begin{smallmatrix}|D|\to 0\\|D'|\to 0\end{smallmatrix}}
\sum_{\begin{smallmatrix}\xi_{i}\in D\\\zeta_{j}\in D'\end{smallmatrix}}
\left[R_{1}\left(\begin{array}{c}\nu_{i}\\1\end{array},\begin{array}{c}a\\b\end{array}\right)-
R_{1}\left(\begin{array}{c}0\\\nu_{i}\end{array},\begin{array}{c}a\\b\end{array}\right)\right]
R_{2}\left(\begin{array}{c}\xi_{i}\\\xi_{i+1}\end{array},\begin{array}{c}\zeta_{j}\\\zeta_{j+1}\end{array}\right)\nonumber\\
&&\hskip-4cm=\frac{1}{2}
\lim_{|D|\to 0}
\sum_{\xi_{i}\in D}
\left[R_{1}\left(\begin{array}{c}\nu_{i}\\1\end{array},\begin{array}{c}a\\b\end{array}\right)-
R_{1}\left(\begin{array}{c}0\\\nu_{i}\end{array},\begin{array}{c}a\\b\end{array}\right)\right]
R_{2}\left(\begin{array}{c}\xi_{i}\\\xi_{i+1}\end{array},\begin{array}{c}c\\d\end{array}\right)\label{discr_int}\ ,
\end{eqnarray}

\noindent where $\nu_{i}\in[\xi_{i},\xi_{i+1}]$. On the other hand, from Lemma \ref{approx_norm} we have that
\begin{eqnarray}
\langle f_{n},g\rangle_{H^{\otimes 2}}
&=&\frac{1}{2}\sum_{\begin{smallmatrix}k,l=0\\k>l\end{smallmatrix}}^{2^{n}-1}
\iiiint_{I_{k}^{n}\times[a,b]\times I_{l}^{n}\times[c,d]}dR_{1}(s,t)dR_{2}(u,v)
-\frac{1}{2}\sum_{\begin{smallmatrix}k,l=0\\k<l\end{smallmatrix}}^{2^{n}-1}
\iiiint_{I_{k}^{n}\times[a,b]\times I_{l}^{n}\times[c,d]}dR_{1}(s,t)dR_{2}(u,v)\nonumber\\
&=&\frac{1}{2}\sum_{\begin{smallmatrix}k,l=0\\k>l\end{smallmatrix}}^{2^{n}-1}
R_{1}\left(\begin{array}{c}t_{k}^{n}\\t_{k+1}^{n}\end{array},\begin{array}{c}a\\b\end{array}\right)
R_{2}\left(\begin{array}{c}t_{l}^{n}\\t_{l+1}^{n}\end{array},\begin{array}{c}c\\d\end{array}\right)
-\frac{1}{2}\sum_{\begin{smallmatrix}k,l=0\\k<l\end{smallmatrix}}^{2^{n}-1}
R_{1}\left(\begin{array}{c}t_{k}^{n}\\t_{k+1}^{n}\end{array},\begin{array}{c}a\\b\end{array}\right)
R_{2}\left(\begin{array}{c}t_{l}^{n}\\t_{l+1}^{n}\end{array},\begin{array}{c}c\\d\end{array}\right)\nonumber\\
&=&\frac{1}{2}\sum_{l=0}^{2^{n}-1}
\left[R_{1}\left(\begin{array}{c}t_{l+1}^{n}\\1\end{array},\begin{array}{c}a\\b\end{array}\right)-
R_{1}\left(\begin{array}{c}0\\t_{l}^{n}\end{array},\begin{array}{c}a\\b\end{array}\right)\right]
R_{2}\left(\begin{array}{c}t_{l}^{n}\\t_{l+1}^{n}\end{array},\begin{array}{c}c\\d\end{array}\right)\label{discr_inner}\ .
\end{eqnarray}

\noindent Note that in equation (\ref{discr_int}) we could replace the first $\nu_{i}$ by $\nu_{i}^{1}$ and the second $\nu_{i}$ by $\nu_{i}^{2}$, where $\nu_{i}^{1},\nu_{i}^{2}\in[\xi_{i},\xi_{i+1}]$. This is because the integral (\ref{int_aux_1}) could be split into two. Indeed you could choose $\nu_{i}^{1}=\xi_{i}$ and $\nu_{i}^{2}=\xi_{i+1}$. If you do this, then equation (\ref{discr_inner}) becomes a particular election of the partition $D$ in equation (\ref{discr_int}) and thus
\begin{equation*}
\lim_{n\to\infty}\langle f_{n},g\rangle_{H^{\otimes 2}}=
\iiiint_{[0,1]^{4}}f_{12}^{\mathcal{L}}(s,u)g_{12}(t,v)dR_{1}(s,t)dR_{2}(u,v)
\end{equation*}

\noindent from where the result follows.\quad $\blacksquare$
\subsection{The case of the fractional Brownian motion}
One case of special interest is to explore what happens with the generalised L{\'e}vy area for a fractional Brownian motion (fBm) with Hurst parameter $H\in(0,1)$. It is known that when fBm has Hurst parameter $H\leq\frac{1}{2}$ then its covariance function is of finite $\frac{1}{2H}$-variation, and when $H>\frac{1}{2}$ then its covariance function is of bounded variation. Then the L{\'e}vy area can be defined for two fBm with the same Hurst parameter $H$ as long as $H>\frac{1}{4}$. Moreover we can let one of the processes be as irregular as desired provided that the other one is regular enough. That is, we can let one of the fBm be of Hurst parameter $H<\frac{1}{2}$ as long as the other independent fBm has Hurst parameter $H'>\frac{1}{2}-H$.

\section{Carleman-Fredholm representation}

Consider an isonormal Gaussian process $X=\{X(h)\ |\ h\in H\}$. If $H=L^{2}(T,\mathcal{B},\mu)$, then $H^{\odot2}$ is the space of symmetric square integrable functions $L_{\text{sym}}^{2}(T^{2})$. To each element $f\in H^{\odot2}$ corresponds a symmetric Hilbert--Schmidt operator $F:H\to H$ defined by
\begin{equation*}
(F(h))(t)=\int_{T}f(s,t)h(s)\mu(ds)\ .
\end{equation*}

\noindent Denote by $\{\alpha_{n}\ |\ n\geq 1\}$ the eigenvalues of $F$ repeated according to its multiplicity. Then
\begin{equation*}
\mathbb{E}[e^{zI_{2}(f)}]=\frac{1}{\sqrt{\prod_{n=1}^{\infty}(1-2z\alpha_{n})e^{2z\alpha_{n}}}}
\qquad \text{ for }\ 2|\Re(z)|\sigma(F)<1\ ,
\end{equation*}

\noindent where $\sigma(F)=\sup_{n\geq 1}\{|\alpha_{n}|\}$. The infinite product of the above equality is called a generalised determinant or a Carleman--Fredholm determinant of $F$ (see Dunford and Schwartz \cite[p. 1036]{DS}). There are many works about the characteristic functions of quadratic Wiener functionals. In our case we are interested in a particular functional viewed in different Wiener spaces, therefore the aim of this section is to explicit the procedure to compute the eigenvalues as much as possible. As an example, we will compute the eigenvalues for the particular case where the stochastic processes $X_{1}$ and $X_{2}$ are of the form
\begin{equation*}
X_{i}=\int_{0}^{1}f(s)dW_{i}(s)\qquad f\in L^{2}([0,1])\ ,
\end{equation*}

\noindent where $W_{1}$ and $W_{2}$ are two independent Brownian motion.

It is worthwile to state the relationship between elements of the second Wiener chaos and Hilbert--Schmidt operators by working out a particular example and then extending the results to the general case. To this end we study the case of the L{\'e}vy area for two standard Wiener processes, as we commented in the introduction. For $R_{1}(s,t)=R_{2}(s,t)=s\wedge t$ the Hilbert space $H$ is isometric to $L^{2}([0,1]\times\{1,2\})$. Then the L{\'e}vy kernel $f^{\mathcal{L}}$ defines the Hilbert-Smith operator
\begin{eqnarray*}
F:L^{2}([0,1]\times\{1,2\})&\to&L^{2}([0,1]\times\{1,2\})\\
h&\mapsto&\int_{[0,1]\times\{1,2\}}f_{ij}^{\mathcal{L}}(s,t)h_{i}(s)ds\otimes\text{Card}\ ,
\end{eqnarray*}

\noindent which is reduced to the form
\begin{equation}\label{int_eq_vaps}
F(h)_{j}(t)=
\frac{\delta_{2j}}{2}\left(\int_{0}^{t}h_{1}(s)ds-\int_{t}^{1}h_{1}(s)ds\right)
-\frac{\delta_{1j}}{2}\left(\int_{0}^{t}h_{2}(s)ds-\int_{t}^{1}h_{2}(s)ds\right)
\ .
\end{equation}

\noindent If $h$ is an eigenvector of eigenvalue $\alpha$, $F(h)=\alpha h$, then
 it is continuous because it is defined by an integral, and applying again
 the same argument it is differentiable. Then we differentiate the above expression and obtain the matrix representation
\begin{equation}\label{matrix_eigen}
\left(\begin{array}{c}h_{1}'(t)\\h_{2}'(t)\end{array}\right)=\frac{1}{\alpha}
\left(\begin{array}{cc}0&-1\\1&0\end{array}\right)
\left(\begin{array}{c}h_{1}(t)\\h_{2}(t)\end{array}\right)
=\frac{1}{\alpha}M
\left(\begin{array}{c}h_{1}(t)\\h_{2}(t)\end{array}\right)
\end{equation}

\noindent with solution given by $h(t)=e^{M\frac{t}{\alpha}}h(0)$. From (\ref{int_eq_vaps}) it is clear that $h(1)+h(0)=0$ and thus the eigenvalues satisfy the equation
\begin{equation*}
e^{M\frac{1}{\alpha}}=
\left(\begin{array}{cc}\cos(\alpha^{-1})&-\sin(\alpha^{-1})\\\sin(\alpha^{-1})&\cos(\alpha^{-1})\end{array}\right)=
\left(\begin{array}{cc}-1&0\\0&-1\end{array}\right)=-Id_{2}\ .
\end{equation*}

\noindent Therefore the eigenvalues, $\alpha_{n}$, are $\{\pm(\pi(2n+1))^{-1}\ |\ n\in\mathbb{N}\}$ with multiplicity $2$ since the space of solutions of the ordinary differential equation has dimension $2$. Finally we compute the Carleman-Fredholm determinant to obtain $\mathbb{E}[e^{itA}]=\cosh(t)^{-1}$.

Now we use the same sort of ideas into the abstract setting presented in
Section \ref{wiener_chaos_sec}. Let $f\in H^{\otimes 2}$ and define the operator
$F:=\Psi\circ\Phi_{f}$, such that
\begin{equation*}
F:H\xrightarrow{\Phi_{f}} H^{*}\xrightarrow{\Psi} H\ ,
\end{equation*}

\noindent where $\Psi$ is the duality isomorphism and for $g\in H$ we define $\Phi_{f}(g):H\to\mathbb{R}$ as
\begin{equation*}
\Phi_{f}(g)(h):=\langle f,g\otimes h\rangle_{H^{\otimes 2}},\qquad h\in H\ .
\end{equation*}

\noindent It can be proved that $F$ is a Hilbert-Smith operator.
Note that $g\in H$ is an eigenvector of the operator $F$ with eigenvalue $\alpha$
if and only if $\langle f,g\otimes h\rangle_{H^{\otimes 2}}=\alpha\langle g, h\rangle_{H}$ for all $h\in H$. From Lemma \ref{lemma_iso} we can identify $g(t,i)=\delta_{1i}g_{1}(t)+\delta_{2i}g_{2}(t)$ where $g_{i}\in H_{i}$ and the same sort of identification is valid for $h$, then $g\in H$ is an eigenvector of eigenvalue $\alpha$ if and only if
\begin{equation*}
\sum_{i,j=1}^{2}
\langle f_{ij},g_{i}\otimes h_{j}\rangle_{H_{i}\otimes H_{j}}=
\alpha\sum_{i=1}^{2}\langle g_{i},h_{i}\rangle_{H_{i}}\quad
\forall h_{1}\in H_{1},\
\forall h_{2}\in H_{2}\ .
\end{equation*}

We will say that two covariance functions $R_{1}$ and $R_{2}$ are equivalent
if the associated Hilbert spaces $H_{1}$ and $H_{2}$ are the same. Under this
symmetry of the processes we recover the spectrum structure of the classical L{\'e}vy area.

\begin{proposition}
Under the notation of the previous sections, let $\{X_{1}(t)\ |\ 0\leq t \leq1\}$ and
$\{X_{2}(t)\ |\ 0\leq t \leq1\}$ be continuous centered independent Gaussian processes with equivalent covariance functions $R_{1}$ and $R_{2}$ respectively. Then the corresponding Hilbert--Schmidt operator has eigenvalues with even multiplicity and symmetric with respect to zero. As a consequence the characteristic function of the generalised L{\'e}vy area is of the form
\begin{equation*}
\varphi(t)=\prod_{n\geq 1}\frac{1}{(1+4\alpha^{2}_{n}t^{2})^{m_{n}}}\ ,
\end{equation*}

\noindent where $m_{n}\geq1$.

\end{proposition}

\noindent{\it Proof.}
From the factorisation (\ref{factorisation_H}) it is clear that the symmetry of the
 approximation of the L{\'e}vy kernel $\{f_{n}\}_{n\ge 1}$ is transferred
 to $\hat{f}^{\mathcal{L}}$.
 Then, from  equation
\begin{equation}\label{eigenvalues}
\sum_{\begin{smallmatrix}i,j=1\\i\neq j\end{smallmatrix}}^{2}
\langle \hat{f}_{ij}^{\mathcal{L}},g_{i}\otimes h_{j}\rangle_{H_{i}\otimes H_{j}}=
\alpha\sum_{i=1}^{2}\langle g_{i},h_{i}\rangle_{H_{i}}\quad
\forall h_{1},h_{2}\in H_{1}\equiv H_{2}\ ,
\end{equation}

\noindent it is  checked that if $g(t,i)=\delta_{1i}g_{1}(t)+\delta_{2i}g_{2}(t)$ is an
 eigenvector with eigenvalue $\alpha$, then $\tilde{g}(t,i)=\delta_{1i}g_{2}(t)-\delta_{2i}g_{1}(t)$
 is an eigenvector with eigenvalue $\alpha$ and
 $\hat{g}(t,i)=\delta_{1i}g_{2}(t)+\delta_{2i}g_{1}(t)$ is an eigenvector with eigenvalue $-\alpha$.
  If $g=\lambda \tilde{g}$ for $\lambda\in\mathbb{R}\setminus\{0\}$ then
  $g_{1}=\lambda g_{2}=-\lambda^{2}g_{1}$ and hence $g\equiv0$,
  thus $g$ and $\tilde{g}$ are linear independent. On the other hand,
  if $\{g_{k}(t,i),\tilde{g}_{k}(t,i)\ |\ k=1,\dots, K\}$
   is a family set of linear independent eigenvectors of eigenvalue $\alpha$, and
\begin{equation*}
h(t,i)=\sum_{k=1}^K \lambda_{k}g_{k}(t,i)+\sum_{k=1}^K\mu_{k}\tilde{g}_{k}(t,i),
\qquad \lambda_{k},\mu_{k}\in\mathbb{R},
\end{equation*}

\noindent then
\begin{equation*}
\tilde{h}(t,i)=\sum_{k=1}^K\lambda_{k}\tilde{g}_{k}(t,i)-\sum_{k=1}^K\mu_{k}g_{k}(t,i)\ .
\end{equation*}

\noindent Therefore $\alpha$ has even multiplicity. Note that this suffices to deduce the same property for the eigenvalue $-\alpha$ and by construction the multiplicity of $\alpha$ and $-\alpha$ is the same. Finally, we recover the same structure for the spectrum of the Hilbert--Schmidt operator that we have in the classical L{\'e}vy area. \quad $\blacksquare$

{\bf Example. }
From the explicit calculations made for the classical L{\'e}vy area we can easily work out a bit more general case. Let $X_{i}(t)=\int_{0}^{t}f(s)dW_{i}(s)$, for two independent Brownian motions $W_{1}$ and $W_{2}$. Then $R_{i}(s,t)=\int_{0}^{s\wedge t}f^{2}(u)du$, $H_{i}=L^{2}([0,1],f^{2}(u)du)$ and equation (\ref{matrix_eigen}) can be written as
\begin{equation*}
\left(\begin{array}{c}h_{1}'(t)\\h_{2}'(t)\end{array}\right)=\frac{f^{2}(t)}{\alpha}
\left(\begin{array}{cc}0&-1\\1&0\end{array}\right)
\left(\begin{array}{c}h_{1}(t)\\h_{2}(t)\end{array}\right)\ .
\end{equation*}

\noindent Therefore the general solution is
\begin{equation*}
h(t)=\exp
\left(\int_{0}^{t}\frac{f^{2}(u)}{\alpha}du
\left(\begin{array}{cc}0&-1\\1&0\end{array}\right)
\right)h(0)\ .
\end{equation*}

\noindent Finally the characteristic function of $A$ in this setting is $\mathbb{E}[e^{itA}]=\text{sech}\left(t||f||_{L^{2}}^{2}\right)$.

\subsection*{Acknowledgements.} The authors are partially supported by grant MTM2009-08869 from the
Ministerio de Ciencia e Innovaci{\'o}n and FEDER. A. Ferreiro-Castilla is also supported by a PhD
grant of the Centre de Recerca Matem{\`a}tica.


\end{document}